\documentclass[12pt, a4paper]{article}
\usepackage[latin1]{inputenc}
\usepackage[english]{babel}
\usepackage{indentfirst}
\usepackage{amstext}
\usepackage{amsfonts}
\usepackage{textcomp}
\usepackage{amssymb}
\usepackage{amscd}
\usepackage{graphicx}
\usepackage{setspace}
\usepackage{psfrag,epsf,epsfig}
\usepackage{amsmath}

\newcommand {\demo}{\hskip -0.6cm{\bf Proof.  }}
\newcommand {\fim}{\hfill{$\square$}\vskip 1pc}

\newcommand {\R}{\mathbb{R}}
\newcommand {\N}{\mathbb{N}}

\newcommand {\E}{\mathbb{E}}

\newcommand {\C}{\mathbb{C}}

\newtheorem{teorema}{Theorem}[section]

\newtheorem{definicao}[teorema]{Definition}
\newtheorem{proposicao}[teorema]{Proposition}
\newtheorem{obs}[teorema]{Remark}

\begin{document}

\onehalfspace

\title{Unitary equivalence of representations of graph algebras and branching systems}
\maketitle
\begin{center}

{\large Daniel Gonçalves and Danilo Royer}\\
\end{center}  
\vspace{8mm}

\abstract In this paper we show that, for a class of countable graphs, every representation of the associated graph algebra in a separable Hilbert space is unitarily equivalent to a representation obtained via branching systems.

\section{Introduction}

Graph C*-algebras have been the focus of constant research in recent years. One of the reasons for such interest is the fact that graph C*-algebras, being generalizations of Cuntz algebras, are related to the theory of Wavelets, see \cite{bratteli}. In fact, it is the representations of the Cuntz algebras that are closely related to the theory of wavelets, and hence the study of representations of graph algebras is of interest. 

Recently, we have showed how to obtain representations of a graph algebra in $L^2(\R)$ via branching systems, see \cite{perrongraph} and \cite{perroncuntz}. In this paper we show that, for a large class of graph C*-algebras, all representations of such algebras are unitarily equivalent to a representation arising from a branching system. This class of C*-algebras includes the compact operators and the algebras associated to Bratteli diagrams (with the direction of the edges inverted), among others. The paper is organized as follows: In section 2, we show that, under a technical condition, all representations of a graph algebra are unitarily equivalent to a representation arising from a branching system. With this in mind, our next goal is to find all graph algebras for which we may apply the results of section 2. In order to make the ideas precise, we need to introduce some new terminology for graphs and this is done in section 3. Finally, in section 4, we give a sufficient condition over the graph, to guarantee that all representations of the associated graph C*-algebra are unitary equivalent to a representation arising from a branching system. Before we proceed, we remind the reader of some key definitions and results below (the reader may find more details in \cite{perrongraph} or \cite{perroncuntz}). 

Let $E=(E^0, E^1, r, s)$ be a directed graph, that is, $E^0$ is a set of vertices, $E^1$ is a set of edges and $r,s:E^1\rightarrow E^0$ are the range and source maps. Following [\ref{flr}], the C*-algebra of the graph $E$ is the universal C*-algebra, $C^*(E)$, generated by projections $\{P_v\}_{v\in E^0}$ and partial isometries $\{S_e\}_{e\in E^1}$ with orthogonal ranges satisfying:
\begin{itemize}
\item the projections $p_v$ are mutually orthogonal,
\item $S_e^*S_e=P_{r(e)}$ for each $e\in E^1$,
\item $S_eS_e^*\leq P_{s(e)}$ for each $e\in E^1$,
\item $P_v=\sum\limits_{e:s(e)=v}S_eS_e^*$ for every vertex $v$ with $0<\#\{e:s(e)=v\}<\infty$.
\end{itemize}

For measurable subsets $A,B$ in a given measure space $(X,\mu)$, the notation $B\stackrel{\mu-a.e.}{\subseteq}A$ means that $\mu(B\setminus A)=0$, and the notation $A\stackrel{\mu-a.e.}{=}B$ means that $\mu(A\setminus B)=0$ and $\mu(B\setminus A)=0$. For two maps, $f,g:A\rightarrow X$, the notation $f\stackrel{\mu-a.e.}{=}g$ means that $\mu(x\in A:f(x)\neq g(x))$=0.

\begin{definicao}\label{brancsystem}
Let $(X,\mu)$ be a measure space and let $\{R_e\}_{e\in E^1}$, $\{D_v\}_{v\in E^0}$ be families of measurable subsets of $X$ such that:
\begin{enumerate}
\item $R_e\cap R_d\stackrel{\mu-a.e.}{=} \emptyset$ for each $d,e\in E^1$ with $d\neq e$,
\item $D_u\cap D_v\stackrel{\mu-a.e.}{=} \emptyset$ for each $u,v\in E^0$ with $u\neq v$,
\item $R_e\stackrel{\mu-a.e.}{\subseteq}D_{s(e)}$ for each $e\in E^1$,
\item $D_v\stackrel{\mu-a.e.}{=} \bigcup\limits_{e:s(e)=v}R_e$\,\,\,\,\, if\,\,\,\,\, $0<\#\{e\in E^1\,\,:\,\,s(e)=v\}< \infty$,
\item for each $e\in E^1$, there exists a map $f_e:D_{r(e)}\rightarrow R_e$ such that $f_e(D_{r(e)})\stackrel{\mu-a.e.}{=}R_e$ and the Radon-Nikodym derivative $\Phi_{f_e}$ of $\mu\circ f_e$, with respect to $\mu$ (in $D_{r(e)}$), exists and $\Phi_{f_e}> 0$ $\mu$ a.e.,
\item for each $f_e$ as above there exists a map $f_e^{-1}:R_e \rightarrow D_{r(e)}$ such that $f_e\circ f_e^{-1}\stackrel{\mu-a.e.}{=}Id_{R_e}$ and  $f_e^{-1}\circ f_e\stackrel{\mu a.e.}{=}Id_{D_{r(e)}}$, and for each such $f_e^{-1}$ there exists the Radon-Nikodym derivative $\Phi_{f_e^{-1}}$ of $\mu\circ f_e^{-1}$ with respect to $\mu$ (in $R_e$).

A measurable space $(X,\mu)$, with families of measurable subsets $\{R_e\}_{e\in E^1}$ and $\{D_v\}_{v\in E^0}$, and maps $f_e$, $f_e^{-1}$, $\Phi_{f_e}$ and  $\Phi_{f_e^{-1}}$ as above is called an $E$-branching system.
\end{enumerate}
\end{definicao}

 \begin{teorema}\label{rep}(\cite{perrongraph}) Let $(X,\mu)$ be an $E$-branching system. Then there exists a *-homomorphism $\pi:C^*(E)\rightarrow B(L^2(X,\mu))$ such that 
$$\pi(S_e)\phi=\chi_{R_e}\cdot\Phi_{f_e^{-1}}^{\frac{1}{2}}\cdot \phi\circ f_e^{-1} \text{ and }\pi(P_v)\phi=\chi_{D_v}\phi,$$ for each $e\in E^1$ and $v\in E^0$.
  \end{teorema}

  From now on we suppose that all the graphs are countable, that is, the set of the vertices, $E^0$,  and the set of the edges, $E^1$, are both countable.

\section{Unitary Equivalence of Representations}

In this section we show that, under a technical condition, all representations of a graph C*-algebra are unitarily equivalent to a representation arising from a branching system. To do this, we first need to write the relations defining a graph C*-algebra in terms of relations between the initial and final space of the partial isometries defining the algebra.

Let $H$ be a separable Hilbert space and $\pi:C^*(E)\rightarrow B(H)$ be a representation of the graph algebra $C^*(E)$. Then, the families $\{\pi(S_e)\pi(S_e)^*\}_{e\in E^1}$, $\{\pi(S_e)^*\pi(S_e)\}_{e\in E^1}$, $\{\pi(P_v)\}_{v\in E^0}$ are families of projections in $B(H)$. For each $e\in E^1$, define $H_e=\pi(S_e)\pi(S_e)^*(H)$ and for each $v\in E^0$, define $H_v=\pi(P_v)(H)$, which are closed subspaces of $H$. From the relations that define $C^*(E)$, we obtain:
\begin{enumerate}
\item $H_e\cap H_f=0$ for each $e,f\in E^1$ such that $e\neq f$;
\item $H_u\cap H_v=0$ for each $u,v\in E^0$ such that $u\neq v$;
\item $\pi(S_e):H_{r(e)}\rightarrow H_e$ is isometric and surjective;
\item if $\#\{s^{-1}(v)\}>0$ then $H_v=\left(\bigoplus\limits_{e:s(e)=v}H_e \right) \bigoplus V_v$, where $V_v$ is a subspace of $H_v$. If $0<\#\{s^{-1}(v)\}<\infty$ then $H_v=\bigoplus\limits_{e:s(e)=v}H_e$.

\item $H=\left(\bigoplus\limits_{v\in E^0}H_v\right)\bigoplus V$, where $V$ is a Hilbert space.

\end{enumerate}  

For each $e\in E^1$, choose an orthonormal (Schauder) basis $B_e$ of $H_e$. Since $H$ is separable, $B_e$ is finite or countable.
For the subspaces $H_v$, with $v\in E^0$, choose a basis as follows:
\begin{itemize}
\item if $\#\{s^{-1}(v)\}=\infty$ then $H_v=\left(\bigoplus\limits_{e:s(e)=v}H_e\right)\bigoplus V_v$, and choose an orthonormal (Schauder) basis $B_v$ of $H_v$ such that if $e\in s^{-1}(v)$ then $B_e\subseteq B_v$. If $0<\#\{s^{-1}(v)\}<\infty$ then choose the basis $B_v$ of $H_v$ as being $B_v:=\bigcup\limits_{e\in s^{-1}(v)}B_e$.

\item if $\#\{s^{-1}(v)\}=0$, choose any orthonormal (Schauder) basis $B_v$ of $H_v$.
\end{itemize}

Finally, choose an orthonormal (Schauder) basis $B$ of $H$ such that $B_v\subseteq B$ for each $v\in E^0$. Notice that such basis exists because $H=\left(\bigoplus\limits_{v\in E^0}H_v\right)\bigoplus V$, and moreover $B$ is countable because $E^0$ and $E^1$ are both countable.

To show unitarily equivalence between a representation of a graph algebra and a representation arising from a branching system we have to assume that the unitary operator, $\pi(S_e): H_{r(e)}\rightarrow H_e$, takes the basis of $H_{r(e)}$ to the basis of $H_e$. Namely, we have to ask that

\begin{equation}\label{h1}\tag{$B2B$}
\pi(S_e)(B_{r(e)})=B_e ,\text{ for each } e\in E^1,
\end{equation}
which we call condition (\ref{h1}). Notice that condition (\ref{h1}) is equivalent to say that $\pi(S_e)^*(B_e)=B_{r(e)}$ for each $e\in E^1$.

We may now prove the following:

\begin{teorema}\label{unitequivalent}
Let $\pi:C^*(E)\rightarrow B(H)$ be a *-representation of $C^*(E)$, where $H$ is separable. Choose an orthonormal basis $B=\{h_j\}_{j\in \N}$ of $H$ as above and suppose this basis satisfies the hypothesis (\ref{h1}). Then the representation $\pi$ is unitarily equivalent to a representation $\tilde{\pi}:C^*(E)\rightarrow B(l^2(\N))$, where $\widetilde{\pi}$ is induced by an $E-$branching system.
\end{teorema}

\demo Let $l^2(\N)=\{(x_n)_{n\in\N}:x_n\in \C\,\, \forall n\in \N \text{ and } \sum\limits_{n\in \N}|x_n|^2<\infty\}$ and let $\{d^n\}_{n\in \N}$ be the canonical basis in $l^2(\N)$. 
Define $U:H\rightarrow l^2(\N)$ by $U(h_j)=d^j$, which is an unitary operator. Next, we will define the $E$-branching system $(X,\mu)$. For this, let $X=\N$ and let $\mu$ be the counting measure in $\N$. For each $e\in E^1$, define  
$$R_e=\{j\in \N: h_j\in B_e\}$$ and for each $v\in E^0$ define 

$$D_v=\{j\in\N:h_j\in B_v\}.$$
To check that $(X, \mu)$, with the families $\{D_v\}_{v\in E^0}$ and $\{R_e\}_{e\in E^1}$, is an $E$-branching system, we need to verify the conditions of $\ref{brancsystem}$. Since $H_e\cap H_f=0$ then $R_e\cap R_f=\emptyset$ for $e\neq f$. Similarly, $D_u\cap D_v=\emptyset$. Conditions $3$ and $4$ follows by the choice of $B=\{h_j\}_{j\in \N}$. So it remains to define, for each $e\in E^1$, a map $f_e:D_{r(e)}\rightarrow R_e$ according to definition $\ref{brancsystem}$.

For a given $e\in E^1$, we define $f_e$ as follows: if $j_0\in D_{r(e)}$ then $h_{j_0}\in H_{r(e)}$ and by hypothesis (\ref{h1}), $\pi(S_e)(h_{j_0})=h_{i_0}$ for some $i_0\in \N$. Note that since $h_{i_0}\in H_e$ then $i_0\in R_e$. Define $f_e(j_0)=i_0$, and so $f_e:D_{r(e)}\rightarrow R_e$ is a bijection. (The map $f_e$ is surjective because $\pi(S_e)(B_{r(e)})=B_e$).

Note that in this case, since $\mu$ is the counting measure, the Radon-Nykodim derivatives are $\Phi_{f_e}(x)=1$, for each $x\in D_{r(e)}$, and $\Phi_{f_e^{-1}}(x)=1$ for each $x\in R_e$. Since $(X,\mu)$ is an $E$-branching system, by theorem \ref{rep} there exists a *-homomorphism $\widetilde{\pi}:C^*(E)\rightarrow B(L^2(X,\mu))=B(l^2(\N))$ such that for each $c=(c_j)_{j\in \N}\in l^2(\N)$,   $\widetilde{\pi}(S_e)c=\left((\widetilde{\pi}(S_e)c)_j\right)_{j\in \N}$, where $$(\widetilde{\pi}(S_e)c)_j=[j\in R_e]c_{f_e^{-1}(j)}.$$ ($[j\in R_e]=1$ if $j\in R_e$ and $[j\in R_e]=0$ if $j\notin R_e$). Also, by theorem \ref{rep} for each $v\in E^0$,
$\widetilde{\pi}(P_v)c=((\widetilde{\pi}(P_v)c)_j)_{j\in \N}$ where $$(\widetilde{\pi}(P_v)c)_j=[j\in D_v]c_j.$$

To show that $\pi$ and $\widetilde{\pi}$ are unitarily equivalent, we need to show that $U^*\widetilde{\pi}(S_e)U=\pi(S_e)$, for each $e\in E^1$, and $U^*\widetilde{\pi}(P_v)U=\pi(P_v)$, for each $v\in E^0$. 

First we will show that $U^*\widetilde{\pi}(S_e)U=\pi(S_e)$, for each $e\in E^1$. For this, fix an element $e\in E^1$ and let $h_k$ be a vector of the basis $B$ of $H$. Then, $$\widetilde{\pi}(S_e)U(h_k)=\widetilde{\pi}(s_e)(d^k)=((\tilde{\pi}(s_e)d^k)_j)_{j\in \N}=([j\in R_e]d^k_{f_e^{-1}(j)})_{j\in \N}.$$ 
Notice that $[j\in R_e]d^k_{f_e^{-1}(j)}=1$ if and only if $j\in R_e$ and $f_e^{-1}(j)=k$, or equivalently, if $k\in D_{r(e)}$ and $f_e(k)=j$. So if $h_k\notin B_{r(e)}$ then $\widetilde{\pi}(S_e)(d^k)=0$ and $\pi(S_e)(h_k)=0$, and hence $\pi(S_e)(h_k)=0=U^*\pi(S_e)U(h_k)$. 
If $h_k\in B_{r(e)}$ then $\pi(S_e)(h_k)=h_{i_0}$ for some $i_0$, by (\ref{h1}), and so $f_e(k)=i_0$. In this case, $[j\in R_e]d^k_{f_e^{-1}(j)}=1$ if and only if $j=i_0$, and hence $\widetilde{\pi}(S_e)d^k=d^{i_0}$. Therefore, in this case, $$U^*\widetilde{\pi}(S_e)U(h_k)=U^*\widetilde{\pi}(S_e)d^k=U^*d^{i_0}=h_{i_0}=\pi(S_e)(h_k).$$

So, for each vector $h_k\in B$, $U^*\widetilde{\pi}(S_e)U(h_k)=\pi(S_e)(h_k)$. This is enough to see that $U^*\widetilde{\pi}(S_e)U=\pi(S_e)$.

It remains to show that $U^*\widetilde{\pi}(P_v)U=\pi(P_v)$ for each $v\in E^0$. Choose an element $v\in E^0$. Then, for a vector $h_k\in B$, $\pi(P_v)(h_k)=[h_k\in B_v]h_k$ and $\widetilde{\pi}(P_v)d^k=[k\in D_v]d^k$. So, $$U^*\widetilde{\pi}(P_v)U(h_k)=U^*\widetilde{\pi}(P_v)(d^k)=[k\in D_v]U^*(d^k)=[k\in D_v]h_k=\pi(P_v)(h_k).$$ It follows that $\pi(P_v)=U^*\widetilde{\pi}(P_v)U$ and hence the representation $\pi$ is unitarily equivalent to $\widetilde{\pi}$ as desired.
\fim

\section{Terminology and results for graphs}

In this section we introduce some terminology and definitions for graphs which we will need in the next section (where we give a list of graph algebras for which we may apply the results of the previous section).

Our first definition is related to non oriented paths in a graph:

\begin{definicao} Let $E=(E^0,E^1,r,s)$ be a graph. We say that:
\begin{enumerate}
\item A vertex $u$ is adjacent to an edge $e$, or $e$ is adjacent to $u$, if $r(e)=u$ or $s(e)=u$.
\item Two vertices $u\neq v$ are adjacent if there exists an edge $e$ adjacent to $u$ and $v$.
\item Two edges $e\neq f$ are adjacent if there exists a vertex $u$ such that $u$ and $e$, and $u$ and $f$ are adjacent. 
\item A path between $u,v\in E^0$ is a pair of sequences ($u_0u_1...u_n;e_1...e_n$) of vertices $u_i$ and edges $e_j$ such that $e_i$ is adjacent to $u_{i-1}$ and $u_i$, for each $i\in \{1,...,n\}$, $u=u_0$, $v=u_n$ and $e_i\neq e_j$ for $i\neq j$.
\item A cycle is a path $(u_0...u_n; e_1...e_n)$ such that $u_0=u_n$.

\item A graph $E$ is $P$-simple if for each $u,v \in r(E^1)\cup s(E^1)$, with $u\neq v$, there exists at most one path between $u$ and $v$, and moreover there does not exist $e\in  E^1$ such that $r(e)=s(e)$.
\end{enumerate}
\end{definicao}

Notice that a graph is $P$-simple if and only if it contains no cycles and $\nexists\,\, e \in E^1$ such that $r(e)=s(e)$.

\begin{definicao}\label{connectedset} Let $E=(E^0,E^1, r,s)$ be a graph. We say that a subset $Z$ of $E^0$ is connected if, for each $u,v\in Z$, there exists a path between $u$ and $v$. We say that the connected subset $Z$ is maximal if the following holds: if $v\in E^0$ is such that there exists a path between $v$ and $u$, where $u\in Z$, then $v\in Z$.
\end{definicao}

Notice that if $Z\subseteq E^0$ is connected then necessarily $Z\subseteq r(E^1)\cup s(E^1)$. Also, if $E$ is $P$-simple then $Z\subseteq E^0$ be connected  means that for each pair of distinct vertices $u,v\in Z$ there exists exactly one path between $u$ and $v$.

For a given graph $E=(E^0,E^1,r,s)$, $E^0$ obviously does not need to be connected, but $E^0$ is a disjoint union of connected maximal subsets with the subset of the ``isolated vertices", that is, the vertices which do not belong to $r(E^1)\cup s(E^1)$. 

To obtain one connected maximal subset $Z_v$, fix some vertex $v\in r(E^1)\cup s(E^1)$ and define $$Z_v:=\{v\}\cup\{u\in E^0: \text{ there is a path between $u$ and $v$} \}.$$ Using Zorn's Lemma we may write $$E^0=\left(\bigcup\limits_{i\in\Delta}^.Z_i\right)\bigcup\limits^.R$$ where each $Z_i$ is a connected maximal subset of $E^0$ and $R=E^0\setminus r(E^1)\cup s(E^1)$ is the set of the isolated vertices. 

Notice that since each $Z_i$ is connected and maximal then it holds that $r^{-1}(Z_i)\cup s^{-1}(Z_i)=r^{-1}(Z_i)=s^{-1}(Z_i)$, and so $(r^{-1}(Z_i)\cup s^{-1}(Z_i),Z_i,r_i,s_i)$ is a subgraph of $E$, where the maps $r_i,s_i$ are the restrictions of $r$ and $s$. 

So, except for the set $R$, the graph $E$ is ``the disjoint union" of the subgraphs $(r^{-1}(Z_i)\cup s^{-1}(Z_i),Z_i,r_i,s_i)$, where each $Z_i$ is connected and maximal. 

Before we give some examples, we need one more definition:

\begin{definicao} A vertex $v\in E^0$ is an extreme vertex of $E$ if $\#\{r^{-1}(v)\cup s^{-1}(v)\}=1$ and if there does not exist an edge $e\in E^1$ such that $r(e)=v=s(e)$. If $v$ is an extreme vertex, then the unique edge adjacent to $v$ is called an extreme edge.
\end{definicao} 

 In other words, a vertex $v$ is an extreme vertex of $E$ if $v$ is adjacent to exactly one other vertex and if there does not exist an edge $e$ with $r(e)=v=s(e)$.
Notice that a graph does not need to have extreme vertices. For example, the following graph 

\centerline{
\setlength{\unitlength}{2cm}
\begin{picture}(4,0.6)
\put(0,0){$\dots$}
\put(0.5,0){\circle*{0.05}}
\put(0.6,0){\vector(1,0){1}}
\put(0.4,0.1){$v_{-1}$}
\put(1.1,0.1){$e_0$}
\put(1.7,0){\circle*{0.05}}
\put(1.8,0){\vector(1,0){1}}
\put(1.6,0.1){$v_0$}
\put(2.2,0.1){$e_1$}
\put(2.9,0){\circle*{0.05}}
\put(3,0){\vector(1,0){0.9}}
\put(2.8,0.1){$v_1$}
\put(3.4,0.1){$e_2$}
\put(4,0){\circle*{0.05}}
\put(3.9,0.1){$v_2$}
\put(4.1,0){\dots}
\end{picture}}
\vspace{0.5pc}
whose graph algebra is the algebra of the compact operators has no extreme vertices.

In the next pages we will describe a procedure to classify the edges and vertices of a graph taking in consideration "how extreme" they are. Such classification will be of importance later.

Fix a graph $E$ and suppose $E$ has at least one extreme vertex. 
Let $X_1$ be the set of extreme vertices of $E$ and $Y_1$ be the set of extreme edges of $E$. The vertices in $X_1$ and the edges in $Y_1$ will be called respectively by level 1 vertices and level 1 edges. Consider now the graph $\E_1=(E^1\setminus Y_1,E^0\setminus X_1, r,s)$ where $r,s$ are, by abuse of notation, respectively the restriction maps $r_{|_{E^1\setminus Y_1}}, s_{|_{E^1\setminus Y_1}}:E^1\setminus Y_1\rightarrow E^0 \setminus X_1$.

Next, let $X_2$, $Y_2$ be, respectively, the sets of extreme vertices (supposing there is at least one such a vertex) and extreme edges of the graph $\E_1$. The elements of $X_2$ and $Y_2$ will be called level 2 vertices and level 2 edges of the graph $E$, respectively.

Now define the graph $\E_2$ as being the graph $\E_2=(E^1\setminus(Y_1\cup Y_2), E^0\setminus (X_1\cup X_2), r,s)$. Again, by abuse of notation, $r$ and $s$ mean the restriction of $r,s$ to $E^1\setminus(Y_1\cup Y_2)$. 

By inductive reasoning, we define the sets $X_n$ and $Y_n$ as being the (nonempty) sets of extreme vertices and extreme edges of the graph $\E_{n-1}$, and call respectively the elements of $X_n$ and $Y_n$ by level $n$ vertices and level $n$ edges of $E$. Notice that the level $n$ vertices and level $n$ edges of the graph $E$ are nothing more than level $1$ vertices and level $1$ edges (or extreme vertices and extreme edges) of the graph $\E_{n-1}$.

If in some of these steps, the graph $\E_m$ has no extreme vertices then there does not exist vertices of level $k$ for any $k>m$. In this case, the level $m$ vertices are the vertices of {\it maximum level}.

For example, let $E$ be the graph as in the next diagram. Below we show the diagrams of the graphs $\E_n$ for those graph.

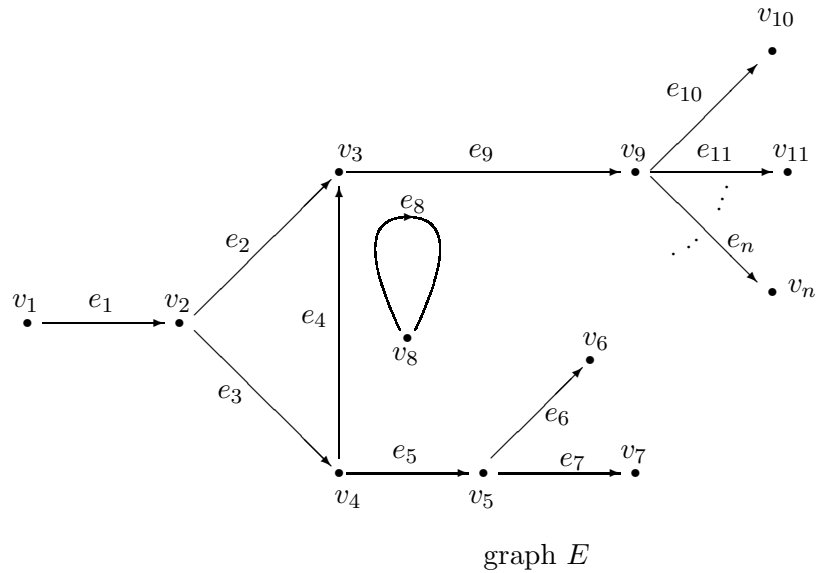
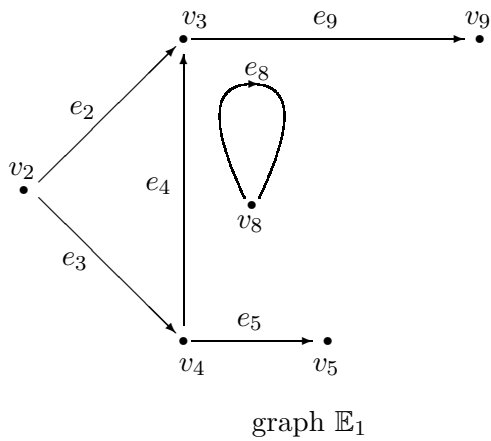
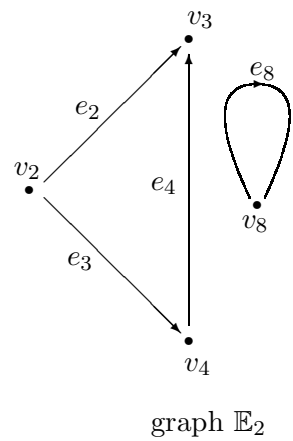
\setlength{\unitlength}{2cm}
\begin{picture}(10,2.3)
\put(0,0){\circle*{0.05}}
\put(0.1,0){\vector(1,0){0.8}}
\put(1,0){\circle*{0.05}}
\put(-0.1,0.1){$v_1$}
\put(0.4,0.1){$e_1$}
\put(0.9,0.1){$v_2$}

\put(2.5,-0.1){\circle*{0.05}}
\put(2.4,-0.25){$v_8$}
\qbezier(2.45,-0.05)(2.1,0.7)(2.5,0.7)
\qbezier(2.55,-0.05)(2.9,0.7)(2.5,0.7)
\put(2.45,0.7){\vector(1,0){0.1}}
\put(2.45,0.76){$e_8$}

\put(1.1,0.05){\vector(1,1){0.9}}
\put(2.05,1){\circle*{0.05}}
\put(1.3,0.5){$e_2$}
\put(2.04,1.1){$v_3$}
\put(1.1,-0.05){\vector(1,-1){0.9}}
\put(2.05,-1){\circle*{0.05}}
\put(1.25,-0.5){$e_3$}
\put(2.02,-1.2){$v_4$}
\put(2.05,-0.9){\vector(0,1){1.8}}
\put(1.8,0){$e_4$}
\put(2.1,-1){\vector(1,0){0.8}}
\put(2.4,-0.9){$e_5$}
\put(3,-1){\circle*{0.05}}
\put(2.9,-1.2){$v_5$}
\put(3.05,-0.9){\vector(1,1){0.6}}
\put(3.4,-0.65){$e_6$}
\put(3.7,-0.25){\circle*{0.05}}
\put(3.65,-0.15){$v_6$}
\put(3.1,-1){\vector(1,0){0.8}}
\put(3.5,-0.95){$e_7$}
\put(4,-1){\circle*{0.05}}
\put(3.9,-0.9){$v_7$}
\put(2.1,1){\vector(1,0){1.8}}
\put(2.9,1.1){$e_9$}
\put(4,1){\circle*{0.05}}
\put(3.9,1.1){$v_9$}
\put(4.1,1){\vector(1,0){0.8}}
\put(4.4,1.1){$e_{11}$}
\put(5,1){\circle*{0.05}}
\put(4.9,1.1){$v_{11}$}
\put(4.1,1.01){\vector(1,1){0.7}}
\put(4.2,1.5){$e_{10}$}
\put(4.9,1.8){\circle*{0.05}}
\put(4.8,2){$v_{10}$}
\put(4.6,0.9){\circle*{0.01}}
\put(4.57,0.82){\circle*{0.01}}
\put(4.55,0.75){\circle*{0.01}}
\put(4.1,0.97){\vector(1,-1){0.7}}
\put(4.6,0.5){$e_n$}
\put(4.9,0.2){\circle*{0.05}}
\put(5,0.2){$v_n$}
\put(4.4,0.55){\circle*{0.01}}
\put(4.35,0.5){\circle*{0.01}}
\put(4.25,0.455){\circle*{0.01}}
\put(3,-1.6){graph $E$}
\end{picture}

\vspace{2.5cm}

\setlength{\unitlength}{2cm}
\centerline{
\begin{picture}(6,2.5)
\put(1.5,-0.1){\circle*{0.05}}
\put(1.4,-0.25){$v_8$}
\qbezier(1.45,-0.05)(1.1,0.7)(1.5,0.7)
\qbezier(1.55,-0.05)(1.9,0.7)(1.5,0.7)
\put(1.45,0.7){\vector(1,0){0.1}}
\put(1.45,0.76){$e_8$}
\put(0,0){\circle*{0.05}}
\put(-0.1,0.1){$v_2$}
\put(0.1,0.05){\vector(1,1){0.9}}
\put(1.05,1){\circle*{0.05}}
\put(0.3,0.5){$e_2$}
\put(1.04,1.1){$v_3$}
\put(0.1,-0.05){\vector(1,-1){0.9}}
\put(1.05,-1){\circle*{0.05}}
\put(0.25,-0.5){$e_3$}
\put(1.02,-1.2){$v_4$}
\put(1.05,-0.9){\vector(0,1){1.8}}
\put(0.8,0){$e_4$}
\put(1.1,-1){\vector(1,0){0.8}}
\put(1.4,-0.9){$e_5$}
\put(2,-1){\circle*{0.05}}
\put(1.9,-1.2){$v_5$}
\put(1.1,1){\vector(1,0){1.8}}
\put(1.9,1.1){$e_9$}
\put(3,1){\circle*{0.05}}
\put(2.9,1.1){$v_9$}
\put(1.5,-1.6){graph $\E_1$}
\put(5.5,-0.1){\circle*{0.05}}
\put(5.4,-0.25){$v_8$}
\qbezier(5.45,-0.05)(5.1,0.7)(5.5,0.7)
\qbezier(5.55,-0.05)(5.9,0.7)(5.5,0.7)
\put(5.45,0.7){\vector(1,0){0.1}}
\put(5.45,0.76){$e_8$}
\put(4,0){\circle*{0.05}}
\put(3.9,0.1){$v_2$}
\put(4.1,0.05){\vector(1,1){0.9}}
\put(5.05,1){\circle*{0.05}}
\put(4.3,0.5){$e_2$}
\put(5.04,1.1){$v_3$}
\put(4.1,-0.05){\vector(1,-1){0.9}}
\put(5.05,-1){\circle*{0.05}}
\put(4.25,-0.5){$e_3$}
\put(5.02,-1.2){$v_4$}
\put(5.05,-0.9){\vector(0,1){1.8}}
\put(4.8,0){$e_4$}
\put(4.8,-1.6){graph $\E_2$}
\end{picture}}

\vspace{3.5cm}

In this example, the graph $\E_2$ has no extreme vertices, and so there are no vertices and edges of level greater or equal to $3$ in the graph $E$. There are only vertices of level $1$ and $2$. The level $1$ vertex set is $X_1=\{v_1,v_6,v_7,v_{10},v_{11},...\}$ (which is infinite), and the maximum level set is $X_2=\{v_5,v_9\}$. The level $1$ edge set is $Y_1=\{e_1,e_6,e_7,e_{10},e_{11},...\}$ and the level $2$ edge set is $Y_2=\{e_5,e_9\}$.

In general, the subsets $X_i$ are nonempty and disjoint. Moreover no vertex not in $r(E^1)\cup s(E^1)$ may be a vertex of some $X_i$, that is, $$\bigcup\limits_{i\in L}X_i\subseteq r(E^1)\cup s(E^1)$$ where $L=\emptyset$ if $E$ has no extreme vertices, $L=\{1,...,m\}$ if $E$ has vertices of maximum level, or $L=\N$ if $E$ has extreme vertices but has no vertices of maximum level.  

Also, in general, it is not true that $$\bigcup\limits_{i\in L}X_i=r(E^1)\cup s(E^1).$$ For example, if $L=\emptyset$ or if there is a cycle in $E$ then the equality above is not verified, as we could see in the previous example. But in some cases the equality above holds, and in these cases some interesting facts may be proved, as we show in the next proposition. 

\begin{proposicao}\label{prop1} Let $E=(E^0,E^1,r,s)$ be a graph, and let $X_i$ be the level $i$ vertex set. Let $Z=s(E^1)\cup r(E^1)$.
\begin{enumerate} 
\item For each level $n$ vertex $v$ there exists at most one vertex $w$ adjacent to $v$ with level greater than or equal to $n$.

\item If $Z=\bigcup\limits_{n=1}^mX_n$ and $Z$ is connected then
\begin{enumerate}
\item if $v\in X_n$ with $n<m$ then there exists exactly one vertex $w$ with level greater than $n$ adjacent to $v$,
\item $X_m$ is a set with two vertices and there exists exactly one edge $e$ adjacent simultaneously to the two vertices.   
\end{enumerate}
\item If $Z=\bigcup\limits_{n=1}^mX_n\bigcup\limits^\cdot \{\overline{v}\}$ and $Z$ is connected then:

\begin{enumerate}
\item if $v\in X_n$ with $n<m$ then there exists exactly one vertex $w$ with level greater than $n$ adjacent to $v$,
\item for each $v\in X_m$ there exists exactly one edge adjacent to $v$ and $\overline{v}$.
\end{enumerate}

\item If $Z$ is finite and connected and if $E$ is P-simple then $Z=\bigcup\limits_{n=1}^mX_n$ or $Z=\left(\bigcup\limits_{n=1}^mX_n\right)\bigcup\limits^\cdot \{\overline{v}\}$.
 \end{enumerate}
\end{proposicao}

\demo 1) Let $v\in X_n$ and suppose $u$ and $w$ are vertices with level greater than or equal to $n$ adjacent to $v$. Let $e_1,e_2$ be edges such that $e_1$ is adjacent to $v$ and $w$ and $e_2$ is adjacent to $v$ and $u$. Then, $e_1$ and $e_2$ are not edges of level less than or equal to $n-1$, and so they belong to the graph $\E_{n-1}$. Therefore, $v$ is not an extreme vertex of $\E_{n-1}$ and this means that $v$ is not a level $n$ vertex.

2) First we prove (a). 

Let $n<m$, $u\in X_n$ and $w\in X_m$. Since $Z$ is connected there exists a path $(u_0...u_p; e_1...e_p)$ (with $e_i\neq e_j$) between $u$ and $w$. If $p=1$ then (a) is proved. So, let $p>1$. Suppose that the level of $u_1$ is less or equal than $n$ and let $u_r$ be (one of) the vertex of smallest level among the vertices $u_1,...,u_p$. Then $u_r$ is adjacent to $u_{r-1}$ and $u_{r+1}$, both of level greater than or equal to the level of $u_r$, which is impossible, by 1. It follows that the level of $u_1$ is greater than $n$.

(b) Let $u,v$ and $w$ be vertices of $X_m$ and let $(u_0...u_p;e_1...e_p)$ be a path from $u$ to $v$. By 1, no vertex, among $u_0,...,u_p$, is a vertex of level less than $m$, (and so each $u_i$ is a level $m$ vertex). So, if $p\geq 2$ then $u_1$ is adjacent to two vertices of level $m$, which is impossible, by 1. Then $p=1$ and so $u$ and $v$ are adjacent. Similarly $v$ and $w$ are adjacent. So, $v$ is adjacent to $u$ and to $w$, which is impossible, again by 1. So, it follows that $X_m$ is a set with one or two vertices. Suppose $X_m=\{u\}$. Then each other vertex $v$ adjacent to $u$ is a vertex of level less than $m$. So, each edge adjacent to $u$ is an edge of level less than $m$. Then there is no edge in $\E_{m-1}$ adjacent to $u$, and this means that $u$ is not an extreme vertex of $\E_{m-1}$, or equivalently, $u\notin X_m$. It follows that $X_m=\{u,v\}$ is a set of two vertices. Since $Z$ is connected there exists a path $(u_0...u_p;e_1...e_p)$ between $u$ and $w$, and by 1, $p=1$, and hence $u$ and $v$ are adjacent. Supposing that there exists two edges adjacent to $u$ and $w$ simultaneously, it follows that $u$ (and $v$) is not an extreme vertex of $\E_{m-1}$, and so $u\notin X_m$, which is a contradiction. So, there exists exactly one edge adjacent to $u$ and $v$.

3) The proof of (a) is the same as the proof of 2(a). 

b) Let $u\in X_m$. Then $u$ is adjacent to exactly one edge $e$ in $\E_{m-1}$. Let $v$ be the other vertex adjacent to $e$. Since $e\in \E_{m-1}$, then either $v\in X_m$ or $v=\overline{v}$. Suppose $v\in X_m$. Since $Z\setminus (X_1\cup...\cup X_{m-1})$ is connected (in $\E_{m-1}$), there exists a path from $v$ to $\overline{v}$ (in $\E_{m-1}$). So, $v$ is adjacent to $u$ and to one more vertex (in $\E_{m-1}$), and hence $v\notin X_m$, which is a contradiction. So $v=\overline{v}$, and hence it follows that $u$ is adjacent to $\overline{v}$. If we suppose that there are two edges adjacent simultaneously to $u$ and to $\overline{v}$, then it follows that $u$ is not an extreme vertex in $\E_{m-1}$, or, equivalently, $u\notin X_m$. So, for each vertex $u\in X_m$ there exists exactly one edge $e$ adjacent to $u$ and to $\overline{v}$.

4) Since $E$ is finite and $P$-simple there exists at least one extreme vertex $v$ of $E$ (supposing $E^1\neq\emptyset$), and hence $X_1\neq\emptyset$. Since $Z$ is finite then there exists some vertex of maximum level. Let $X_m$ be the maximum level set. 

Suppose $Z\neq\bigcup\limits_{n=1}^mX_n$. Since $E$ is P-simple and $Z$ is connected then $\E_m$ is P-simple and $Z\setminus(X_1\cup...\cup X_m)$ is connected (in $\E_m$). So, since $\E_m$ has no extreme vertices (since $X_m$ is the maximum level set), each vertex of $Z\setminus\bigcup\limits_{n=1}^mX_n$ is adjacent to two edges in $\E_m$. Since $Z$ is finite, if we suppose that $X\setminus\bigcup\limits_{n=1}^mX_n$ contains more than one element, we obtain a cycle. But such path does not exist, because $E$ is P-simple. Then it follows that $X\setminus\bigcup\limits_{n=1}^mX_n$ is a set of one element, $\overline{v}$. 
\fim

\begin{obs}
Following proposition \ref{prop1}, if $E$ is a P-simple, connected graph such that $r(E^1)\cup s(E^1)$ is finite, then $$r(E^1)\cup s(E^1)=\bigcup\limits_{n=1}^m X_n$$ or  $$r(E^1)\cup s(E^1)=\left(\bigcup\limits_{n=1}^m X_n\right)\cup\{\overline{v}\}.$$ 
\end{obs}
However, $r(E^1)\cup s(E^1)$ does not need to be finite to have this propriety. For example,
consider the following graph:

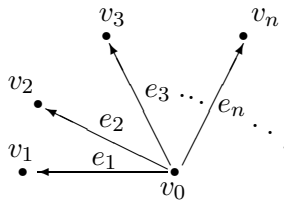
\setlength{\unitlength}{2cm}
\centerline{
\begin{picture}(0,1.5)
\put(0,0){\circle*{0.05}}
\put(-0.05,0){\vector(-1,0){0.85}}
\put(-1,0){\circle*{0.05}}
\put(-1.1,0.1){$v_1$}
\put(-0.1,-0.15){$v_0$}
\put(-0.55,0.05){$e_1$}
\put(-0.05,0.02){\vector(-2,1){0.8}}
\put(-0.9,0.45){\circle*{0.05}}
\put(-1.1,0.55){$v_2$}
\put(-0.5,0.3){$e_2$}
\put(-0.03,0.05){\vector(-1,2){0.4}}
\put(-0.45,0.9){\circle*{0.05}}
\put(-0.5,1){$v_3$}
\put(-0.2,0.5){$e_3$}
\put(0.02,0.5){.}
\put(0.08,0.48){.}
\put(0.15,0.45){.}
\put(0.05,0.05){\vector(1,2){0.4}}
\put(0.45,0.9){\circle*{0.05}}
\put(0.5,1){$v_n$}
\put(0.28,0.4){$e_n$}
\put(0.5,0.3){.}
\put(0.6,0.25){.}
\put(0.7,0.15){.}
\end{picture}}
\vspace{1cm}
This is a connected and P-simple graph, but $r(E^1)\cup s(E^1)$ (which equals to $E^0$ in this case) is not finite. In this case, $X_1=\{v_i:i\geq 1\}$, which is a infinite set, and $r(E^1)\cup s(E^1)=X_1\cup\{v_0\}$.

\section{Graph Algebras whose representations are equivalent to representations arising from Branching Systems.}

According to theorem \ref{unitequivalent}, we can guarantee that a representation of a graph algebra, in a separable Hilbert space $H$, is unitarily equivalent to a representation induced by a branching system if there exists basis of certain subspaces of $H$ (the subspaces $H_e$ and $H_v$) with some particular properties. Next we prove that, under some hypothesis over the graph $E$, there always exists basis as required in theorem \ref{unitequivalent}. Throughout this section we assume the terminology of the previous section.

\begin{teorema}\label{theorem1} Let $E=(E^0,E^1,r,s)$ be a graph such that $Z=r(E^1)\cup s(E^1)$ is connected and suppose $Z=\bigcup\limits_{n=1}^m X_n$ or $Z=\bigcup\limits_{n=1}^m X_n\cup\{\overline{v}\}$, where $X_n$ are as in proposition \ref{prop1}. Let $\pi:C^*(E)\rightarrow B(H)$ be a *-representation, where $H$ is a separable Hilbert space. For each $e\in E^1$ and $v\in E^0$, consider the subspaces $H_e:=\pi(S_eS_e^*)(H)$ and $H_v:=\pi(P_v)(H)$. Then, there exists basis $B_e$ of $H_e$ and $B_v$ of $H_v$ such that:

\begin{enumerate}
\item[1)] if $e\in s^{-1}(v)$ then $B_e\subseteq B_v$ and if $0<|s^{-1}(v)|<\infty$ then $B_v=\bigcup\limits_{e\in s^{-1}(v)}B_e$;
\item[2)] if $e\in r^{-1}(v)$ then $\pi(S_e)(B_v)=B_e$. (and so the basis satisfies hypothesis (\ref{h1})).
\end{enumerate}
\end{teorema}

\demo 
Before we begin the proof of the theorem, let us make some remarks. By proposition \ref{prop1} each $v\in X_n$, with $n<m$, is adjacent to exactly one vertex $u$ with level greater than $n$. If the (unique) edge $e$ adjacent to $v\in X_n$ and to $u$ is such that $r(e)=v$ then we say that $v$ is a {\it final vertex of $X_n$}, and if $s(e)=v$ then we say that $v$ is a {\it initial vertex of $X_n$}. If $X=\bigcup\limits_{n=1}^mX_n$ then by proposition \ref{prop1} $X_m=\{u,v\}$ and there exists exactly one edge $e$ adjacent to $u$ and $v$ simultaneously. Then, $r(e)$ is the {\it final vertex} and $s(e)$ is the {\it initial vertex} of $X_m$. If $X=\bigcup\limits_{n=1}^mX_n\cup \{\overline{v}\}$ then, by proposition \ref{prop1}, for each $v\in X_m$, there exists an unique edge $e$ adjacent to $v$ and $\overline{v}$ simultaneously. In this case, if $r(e)=v$ then $v$ is a {\it final vertex of $X_m$} and if $s(e)=v$ then $v$ is an {\it initial vertex of $X_m$.}

The proof of the theorem will be separated in two steps (Two induction arguments over the level of the vertices). In the first step we will show that there exists basis satisfying conditions 1) and 2) for all final vertices. In the second step we will show how to modify the basis obtained in the first step so that condition 1) and 2) are satisfied for all vertices.

\underline{Step 1:}

Let $v$ be a vertex of level 1 and let $e$ be the (unique) edge adjacent to $v$. If $v$ is a final vertex of $X_1$, then choose a basis $B_v$ of $H_v$ and define a basis $B_e$ of $H_e$ by $B_e:=\pi(S_e)(B_v)$. For all the other vertices $v$ and edges $e$ (of any level), chose any basis $B_v$ and $B_e$. So, for the basis $\{B_v:v\in E^0\}$ and $\{B_e:e\in E^1\}$ the following holds: for each final vertex $v\in X_1$, conditions 1) and 2) are satisfied. Notice that condition 1) is vacuously satisfied.

To proceed with the induction argument, (over $N$), let $N\leq m-1$ and suppose there exists basis $\{\widetilde{B_v}:v\in E^0\}$ and $\{\widetilde{B_e}:e\in E^1\}$, such that for all final vertices $v\in X_k$, with $k\leq N$, conditions 1) and 2) are satisfied.  

To complete the proof we have to show that there exists basis $\{\overline{B_v}:v\in E^0\}$ and $\{\overline{B_e}:e\in E^1\}$ such that properties 1) and 2) holds for each final vertex $v\in X_k$, with $k\leq N+1$. 

So, for each final vertex $v\in X_{N+1}$ proceed as follows: if $s^{-1}(v)=\emptyset$ define $B_v:=\widetilde{B_v}$; if $0<|s^{-1}(v)|<\infty$, define $B_v:=\bigcup\limits_{e\in s^{-1}(v)}\widetilde{B_e}$, and if $|s^{-1}(v)|=\infty$ choose a basis $B_v$ of $H_v$ such that $\widetilde{B_e}\subseteq B_v$ for each $e\in s^{-1}(v)$. In this way we obtain basis $B_v$ for each final vertex $v$ of $X_{N+1}$. Also, for each final vertex $v$ of $X_{N+1}$ and for each $e\in r^{-1}(v)$, define $B_e$ by $B_e:=\pi(S_e)(B_v)$. Finally, for all the other vertices and edges (of any level) define $B_v:=\widetilde{B_v}$ and $B_e:=\widetilde{B_e}$. 

So, the basis $\{B_v:v\in E^0\}$ and $\{B_e:e\in E^1\}$ have the property that if $v\in X_k$, with $k\leq N+1$, is a final vertex then conditions 1) and 2) are satisfied. Hence, by inductive argument, there exists basis $\{\overline{B_v}:v\in E^0\}$ and $\{\overline{B_e}:e\in E^1\}$ such that if $v$ is a final vertex of any level then 1) and 2) are satisfied.

\underline{Step 2:}

Next we will show that the basis above can be modified to a basis satisfying conditions 1) and 2). To do this first we modify the above basis for the initial vertices of $X_m$, and then we show how to redefine the basis for initial vertices of lower levels (via induction).

So, let $\mu$ be the initial vertex of $X_m$ (if $Z=\bigcup\limits_{n=1}^mX_n$) or $\mu=\overline{v}$ (if $Z=\bigcup\limits_{n=1}^mX_n\cup\{\overline{v}\})$. The basis $B_\mu$ will be defined as follows: if $0<|s^{-1}(\mu)|<\infty$, define $B_{\mu}:=\bigcup\limits_{f\in s^{-1}(\mu)}\overline{B_f}$; if $|s^{-1}(\mu)|=\infty$, let $B_\mu$ such that $\overline{B_f}\subseteq B_\mu$ for $f\in s^{-1}(\mu)$ and if $s^{-1}(\mu)=\emptyset$, let $B_\mu:=\overline{B_\mu}$.    

For each $e\in r^{-1}(\mu)$ define $B_e:=\pi(S_e)(B_\mu)$.

For each initial vertex $w$ of $X_m$ (here we are supposing $\mu=\overline{v}$), let $f_0$ be the edge adjacent to $w$ and to $\overline{v}$. Note that $f_0\in s^{-1}(w)$, and $f_0\in r^{-1}(\overline{v})$, since $w$ is an initial vertex of $X_m$. Then, if $|s^{-1}(w)|<\infty$ define $B_w:=B_{f_0}\bigcup\limits_{e\in s^{-1}(w)\setminus\{f_0\}}\overline{B_e}$ and if $|s^{-1}(w)|=\infty$ choose $B_w$ such that $B_{f_0}, \overline{B_e}\subseteq B_w$, for each $e\in s^{-1}(w)\setminus\{f_0\}$.

Finally, for all the other edges $e$ and vertices $v$ define $B_e:=\overline{B_e}$ and $B_v:=\overline{B_v}$. 

So, the basis $\{B_v:v\in E^0\}$, $\{B_e:e\in E^1\}$ are such that for each final vertex $v$ of any level, 1) and 2) are satisfied, and moreover, for each initial vertex $w$ of $X_m$ (and also for $\overline{v}$), 1) and 2) are satisfied. In particular 1) and 2) are satisfied for every vertex of $X_m$ (and hence we proved the first step in the second induction argument).

Now, suppose that there exists basis $\{\widetilde{B_v}:v\in E^0\}$ and $\{\widetilde{B_e}:v\in E^1\}$ such that 1) and 2) are satisfied for each final vertex (of any level) and for all the vertices of level greater or equal to $m-N$ (including the vertex $\overline{v}$), where $m-N\geq 2$.

We need to show that there exists basis $\{B_v:v\in E^0\}$ and $\{B_e:v\in E^1\}$ such that 1) and 2) are satisfied for each final vertex (of any level) and for all the vertices of level greater or equal to $m-(N+1)$. To do this we need to define appropriate basis $B_v$ and $B_e$ for the initial vertices $v$ of $X_{m-(N+1)}$ and for $e\in r^{-1}(X_{m-(N+1)})$.

Let $v$ be an initial vertex of $X_{m-(N+1)}$. If $|s^{-1}(v)|<\infty$, let $B_v:=\bigcup \limits_{e\in s^{-1}(v)}\widetilde{B_v}$  and if $|s^{-1}(v)|=\infty$, define $B_v$ such that $\widetilde{B_f}\subseteq B_v$, for each $f\in s^{-1}(v)$. For each $e\in r^{-1}(v)$, define $B_e:=\pi(S_e)(B_{r(e)})$. For all the other vertices and edges (of any level) define $B_u:=\widetilde{B_u}$ and $B_f:=\widetilde{B_f}$. 

So, $\{B_v:v\in E^0\}$ and $\{B_e:e\in E^1\}$ are basis such that 1) and 2) are satisfied for each final vertex (of any level)  and for all the vertices of level greater than or equal to $m-(N+1)$ and hence we have proved that there exists basis $\{B_v:v\in E^0\}$ and $\{B_e:e\in E^1\}$ satisfying  1) and 2), as desired. 
\fim

For a given graph $E$, according to the remark following definition \ref{connectedset}, $E^0$ may be written as a disjoint union $$E^0=\left(\bigcup^.\limits_{i\in \Delta}Z_i\right)\bigcup\limits^.R,$$ where $\Delta$ is an index set, $R$ is the set of isolated vertices and each $Z_i$ is connected. With this in mind we have the following:

\begin{teorema}\label{thetheorem} 
Let $E=(E^0, E^1,r,s)$ be a graph and write 
$$E^0=\left(\bigcup^.\limits_{i\in \Delta}Z_i\right)\bigcup\limits^.R.$$ Suppose $Z_i=\bigcup\limits_{n=1}^{m_i} X_n$ or $Z_i=\bigcup\limits_{n=1}^{m_i} X_n\cup\{\overline{v_i}\}$ for each $i\in \Delta$, where $X_n$ are as in proposition \ref{prop1}. Then each representation $\pi:C^*(E)\rightarrow B(H)$, where $H$ is a separable Hilbert space, is unitarily equivalent to a representation induced by an $E-$branching system.
\end{teorema} 

\demo The proof follows by applying theorem \ref{theorem1} to each graph $(s^{-1}(Z_i)\cup r^{-1}(Z_i),Z_i,r_i,s_i)$, where $s_i,r_i$ are the restriction maps $r,s:s^{-1}(Z_i)\cup r^{-1}(Z_i)\rightarrow Z_i$, and by theorem \ref{unitequivalent}.
\fim

Notice that the class of graph C*-algebras for which we may define basis satisfying the conditions of theorem \ref{theorem1} is much more extensive then what is covered in theorem \ref{theorem1}. For example, for the algebra of compact operators mentioned before, fixed a vertex $v_0=r(e)$ and chosen a basis for $H_{v_0}$, we can use the isometric isomorphism $\pi(S_e):H_{r(e)} \rightarrow H_e$ to pull back the basis of $H_{v_0}$ to a basis of $H_e$. Then, since $H_v$ should be equal to $\bigoplus\limits_{e:s(e)=v}H_e$, we should define the basis of $H_{s(e)}$ as the basis of $H_e$. Proceeding analogously for $s(e)$ we may "walk to the left" and so on. "Walking to the right" is similar, that is, if $f$ is such that $s(f)=v_0$, then $H_f$ should be equal to $H_{v_0}$ and then $H_{r(f)}$ should be defined using $\pi(S_f)^*$. Following this way one can define all the necessary basis.

Also, for any Bratteli diagram, with the orientation of the edges inverted (so that the first vertex of the diagram is a sink), we may use the techniques presented above to define basis satisfying the conditions of theorem \ref{theorem1}.  We refrain to present a complete argument here, but we believe the reader should be able to apply the ideas presented above to this case as well.

Finally, we leave open the question whether theorem \ref{thetheorem} is valid in general. That is, is every representation of a graph algebra unitary equivalent to a representation arising from a branching system (what branching system would that be?)?

\end{document}